\documentclass{amsproc}
\usepackage{amssymb}
\usepackage{amsmath}
\usepackage{eurosym}
\usepackage{amsfonts}

            \advance \topmargin by -\headheight \advance
            \topmargin by -\headsep
            \evensidemargin
            \oddsidemargin
\newtheorem{theorem}{Theorem}

\newtheorem{remark}[theorem]{Remark}

\numberwithin{equation}{section}
\input{tcilatex}

\begin{document}
\title[Simson Identity of Generalized $m$-step Fibonacci Numbers]{Simson
Identity of Generalized $m$-step Fibonacci Numbers}
\thanks{}
\author[Y\"{u}ksel Soykan]{Y\"{u}ksel Soykan}
\maketitle

\begin{center}
\textsl{Department of Mathematics, Art and Science Faculty, }

\textsl{Zonguldak B\"{u}lent Ecevit University, 67100, Zonguldak, Turkey}

\textsl{e-mail: \ yuksel\_soykan@hotmail.com}
\end{center}

\textbf{Abstract.} One of the best known and oldest identities for the
Fibonacci sequence $\{F_{n}\}$ is 
\begin{equation*}
F_{n+1}F_{n-1}-F_{n}^{2}=(-1)^{n}
\end{equation*}%
which was derived first by R. Simson in 1753 and it is now called as Simson
or Cassini Identity. In this paper, we generalize this result to generalized 
$m$-step Fibonacci numbers and give an attractive formula. Furthermore, we
present some Simson's identities of particular\ generalized $m$-step
Fibonacci sequences.



\section{Introduction}

Several generalizations of Fibonacci numbers and identities have been
studied by mathematicians over the years. In this paper, we generalize
Simson's identity to generalized $m$-step Fibonacci sequences. Before
presenting our main result (Theorem \ref{theorem:hgfdsyunbvgh}) we give some
background. For $m\geq 2,$ the generalized $m$-step Fibonacci numbers, $%
\{V_{n}(V_{0},V_{1},V_{2},...,V_{m-1};r_{1},r_{2},,...,r_{m})\}_{n\geq m}$
(or shortly $\{V_{n}\}_{n\geq m}$), $(n\geq m)$, is defined by the $m$-order
linear recurrence relation%
\begin{equation}
V_{n}=%
\sum_{i=1}^{m}r_{i}V_{n-i}=r_{1}V_{n-1}+r_{2}V_{n-2}+r_{3}V_{n-3}+...+r_{m-1}V_{n-m-1}+r_{m}V_{n-m}
\label{equation:nvxersdxcpomn}
\end{equation}%
with $m$ initial terms%
\begin{equation*}
V_{0}=c_{0},\text{ }V_{1}=c_{1},\text{ }V_{2}=c_{2},...,V_{m-1}=c_{m-1},%
\text{ }
\end{equation*}%
where $r_{i},$ $1\leq i\leq m$, are all real numbers and $c_{i},$ $0\leq
i\leq m-1$, are all real or complex numbers. Such a sequence is also called
the generalized Fibonacci $m$-sequence, or generalized $m$-nacci sequence,
or the $m$-generalized Fibonacci sequence.

The sequences $\{V_{n}\}_{n\geq m}$ can be extended to negative subscripts
by defining 
\begin{equation*}
V_{-n}=-\frac{r_{m-1}}{r_{m}}V_{-(n-1)}-\frac{r_{m-2}}{r_{m}}V_{-(n-2)}-%
\frac{r_{m-3}}{r_{m}}V_{-(n-3)}-...-\frac{r_{1}}{r_{m}}V_{-(n-(m-1))}+\frac{1%
}{r_{m}}V_{-(n-m)}
\end{equation*}%
for $n=m-2,m-1,m,m+1...$. Therefore, recurrence (\ref{equation:nvxersdxcpomn}%
) holds for all integer $n.$

For $m\geq 2,$ the $m$-step Fibonacci numbers, $U_{n}$ $(n\geq m)$, is
defined by the $m$-order linear recurrence relation%
\begin{equation}
U_{n}=\sum_{i=1}^{m}U_{n-i}=U_{m-1}+U_{m-2}+U_{m-3}+...+U_{n-m}
\label{equation:hgftysvfxnb}
\end{equation}%
with $m$ initial terms%
\begin{equation}
\left\{ 
\begin{array}{ccc}
U_{k}=0 & , & -m+2\leq k\leq 0 \\ 
U_{-k+1}=1 & , & k=m%
\end{array}%
\right. .  \label{equation:weybnvfcdxsz}
\end{equation}%
Some of the well known members of this $m$-step Fibonacci numbers include
Fibonacci numbers $F_{n}$ ($m=2,$ $U=F$), Tribonacci numbers $T_{n}$ ($m=3,$ 
$U=T$), Tetranacci numbers $M_{n}$ ($m=4,$ $U=M$) and Pentanacci numbers $%
P_{n}$ ($m=5,$ $U=P$). Here $r_{i}=1$ for all $1\leq i\leq m.$ See Table 1
for some values of these numbers.

\ \ \ \ \ \ 

Table 1 The first few sequences of $m$-step Fibonacci numbers.

\begin{tabular}{cccccccccccccccccccc}
\hline
$m$ & Name & $n$ & $-6$ & $-5$ & $-4$ & $-3$ & $-2$ & $-1$ & $0$ & $1$ & $2$
& $3$ & $4$ & $5$ & $6$ & $7$ & $8$ & $9$ & $10$ \\ \hline
$2$ & Fibonacci & $F_{n}$ & $-8$ & $5$ & $-3$ & $2$ & $-1$ & $1$ & $0$ & $1$
& $1$ & $2$ & $3$ & $5$ & $8$ & $13$ & $21$ & $34$ & $55$ \\ 
$3$ & Tribonacci & $T_{n}$ & $-3$ & $2$ & $0$ & $-1$ & $1$ & $0$ & $0$ & $1$
& $1$ & $2$ & $4$ & $7$ & $13$ & $24$ & $44$ & $81$ & $149$ \\ 
$4$ & Tetranacci & $M_{n}$ & $0$ & $0$ & $-1$ & $1$ & $0$ & $0$ & $0$ & $1$
& $1$ & $2$ & $4$ & $8$ & $15$ & $29$ & $56$ & $108$ & $208$ \\ 
$5$ & Pentanacci & $P_{n}$ & $0$ & $-1$ & $1$ & $0$ & $0$ & $0$ & $0$ & $1$
& $1$ & $2$ & $4$ & $8$ & $16$ & $31$ & $61$ & $120$ & $236$ \\ \hline
\end{tabular}

\ \ \ 

Like the $m$-step Fibonacci numbers, $m$-step Lucas numbers are defined by
the same the $m$-order recurrence relations (\ref{equation:hgftysvfxnb}) but
with different initial terms, namely the $m$-step Lucas numbers, $W_{n}$ ,
is defined by the $m$-order linear recurrence relation%
\begin{equation}
W_{n}=\sum_{i=1}^{m}W_{n-i}=W_{m-1}+W_{m-2}+W_{m-3}+...+W_{n-m}
\end{equation}%
with the $m$ initial terms%
\begin{equation}
\left\{ 
\begin{array}{ccc}
W_{k}=-1 & , & -m+1\leq k\leq -1 \\ 
W_{k}=m & , & k=0%
\end{array}%
\right. .
\end{equation}%
\ \ \ \ 

Some of the well known members of this $m$-step Fibonacci numbers include
Lucas numbers $L_{n}$ ($m=2,$ $W=L$), Tribonacci-Lucas numbers $K_{n}$ ($%
m=3, $ $W=K$), Tetranacci-Lucas numbers $R_{n}$ ($m=4,$ $V=n$) and
Pentanacci-Lucas numbers $Q_{n}$ ($m=5,$ $W=Q$). Here $r_{i}=1$ for all $%
1\leq i\leq m.$ See Table 2 for some values of these numbers.

\ \ \ 

Table 2 The first few sequences of $m$-step Lucas numbers.

\begin{tabular}{cccccccccccccccccc}
\hline
$m$ & Name & $n$ & $-4$ & $-3$ & $-2$ & $-1$ & $0$ & $1$ & $2$ & $3$ & $4$ & 
$5$ & $6$ & $7$ & $8$ & $9$ & $10$ \\ \hline
$2$ & Lucas & $L_{n}$ & $7$ & $-4$ & $3$ & $-1$ & $2$ & $1$ & $3$ & $4$ & $7$
& $11$ & $18$ & $29$ & $47$ & $76$ & $123$ \\ 
$3$ & Tribonacci-Lucas & $K_{n}$ & $-5$ & $5$ & $-1$ & $-1$ & $3$ & $1$ & $3$
& $7$ & $11$ & $21$ & $39$ & $71$ & $131$ & $241$ & $443$ \\ 
$4$ & Tetranacci-Lucas & $R_{n}$ & $7$ & $-1$ & $-1$ & $-1$ & $4$ & $1$ & $3$
& $7$ & $15$ & $26$ & $51$ & $99$ & $191$ & $367$ & $708$ \\ 
$5$ & Pentanacci-Lucas & $Q_{n}$ & $-1$ & $-1$ & $-1$ & $-1$ & $5$ & $1$ & $%
3 $ & $7$ & $15$ & $31$ & $57$ & $113$ & $223$ & $439$ & $863$ \\ \hline
\end{tabular}

\ \ \ \ 

Next we consider the case $r_{i}=1$ for all $1\leq i\leq m-1$ and $r_{m}=2$.
For $m\geq 2,$ $m$-step (order) Jacobsthal numbers, $%
\{J_{n}^{(m)}(J_{0}^{(m)},J_{1}^{(m)},J_{2}^{(m)},...,J_{m-1}^{(m)};1,1,,...,1,2)\}_{n\geq m} 
$ (or shortly $\{J_{n}^{(m)}\}_{n\geq m}$), $(n\geq m)$, is defined by the $%
m $-order linear recurrence relation%
\begin{equation}
J_{n}^{(m)}=\sum_{i=1}^{m-1}r_{i}J_{n-i}^{(m)}+2J_{n-m}^{(m)}
\label{equati:mnvxdzuoaesaoa}
\end{equation}%
with $m$ initial terms 
\begin{equation*}
J_{0}^{(m)}=0\text{ and }J_{i}^{(m)}=1\text{ for }i=1,2,...,m-1.\text{ }
\end{equation*}%
For the $m$th order Jacobsthal-Lucas numbers $j_{n}^{(m)}$ we use the same
recursion (\ref{equati:mnvxdzuoaesaoa}) with initial conditions $%
j_{i}^{(m)}=j_{i}^{(m-1)}$ for $i=0,1,2,...,m-1$ and $j_{0}^{(2)}=2,$ $%
j_{1}^{(2)}=1.$ See Table 3 and Table 4 for $m$th order Jacobsthal numbers
and $m$th order Jacobsthal-Lucas numbers, respectively.

\ \ \ 

Table 3 The first few sequences of $m$th order Jacobsthal numbers.

\begin{tabular}{ccccccccccccccccc}
\hline
$m$ & Name & $n$ & $-3$ & $-2$ & $-1$ & $0$ & $1$ & $2$ & $3$ & $4$ & $5$ & $%
6$ & $7$ & $8$ & $9$ & $10$ \\ \hline
$2$ & second order Jacobsthal & $J_{n}^{(2)}$ & $\frac{3}{8}$ & $-\frac{1}{4}
$ & $\frac{1}{2}$ & $0$ & $1$ & $1$ & $3$ & $5$ & $11$ & $21$ & $43$ & $85$
& $171$ & $341$ \\ 
$3$ & third order Jacobsthal & $J_{n}^{(3)}$ & $-\frac{1}{4}$ & $\frac{1}{2}$
& $0$ & $0$ & $1$ & $1$ & $2$ & $5$ & $9$ & $18$ & $37$ & $73$ & $146$ & $%
293 $ \\ 
$4$ & fourth order Jacobsthal & $J_{n}^{(4)}$ & $\frac{5}{8}$ & $\frac{1}{4}$
& $-\frac{1}{2}$ & $0$ & $1$ & $1$ & $1$ & $3$ & $7$ & $13$ & $25$ & $51$ & $%
103$ & $205$ \\ 
$5$ & fifth order Jacobsthal & $J_{n}^{(5)}$ & $\frac{1}{2}$ & $0$ & $-1$ & $%
0$ & $1$ & $1$ & $1$ & $1$ & $4$ & $9$ & $17$ & $33$ & $65$ & $132$ \\ \hline
\end{tabular}

\ \ \ \ \ \ \ \ 

Table 4 The first few sequences of $m$th order Jacobsthal-Lucas numbers.

\begin{tabular}{ccccccccccccccc}
\hline
$m$ & Name & $n$ & $-3$ & $-2$ & $-1$ & $0$ & $1$ & $2$ & $3$ & $4$ & $5$ & $%
6$ & $7$ & $8$ \\ \hline
$2$ & second order Jacobsthal-Lucas & $j_{n}^{(2)}$ & $-\frac{7}{8}$ & $%
\frac{5}{4}$ & $-\frac{1}{2}$ & $2$ & $1$ & $5$ & $7$ & $17$ & $31$ & $65$ & 
$127$ & $257$ \\ 
$3$ & third order Jacobsthal-Lucas & $j_{n}^{(3)}$ & $1$ & $-1$ & $1$ & $2$
& $1$ & $5$ & $10$ & $17$ & $37$ & $74$ & $145$ & $293$ \\ 
$4$ & fourth order Jacobsthal-Lucas & $j_{n}^{(4)}$ & $-\frac{5}{4}$ & $%
\frac{1}{2}$ & $1$ & $2$ & $1$ & $5$ & $10$ & $20$ & $37$ & $77$ & $154$ & $%
308$ \\ 
$5$ & fifth order Jacobsthal-Lucas & $j_{n}^{(5)}$ & $\frac{1}{4}$ & $\frac{1%
}{2}$ & $1$ & $2$ & $1$ & $5$ & $10$ & $20$ & $40$ & $77$ & $157$ & $314$ \\ 
\hline
\end{tabular}

\ 

For more details about generalized $n$-step Fibonacci numbers we refer to,
for example, the works in [\ref{adegoke2018genfibnum},\ref{bacani2015},\ref%
{cook2013someiden}], among others. Now, we consider the cases $m=2,3,4,5$ of
the generalized $m$-step Fibonacci numbers separately.

Horadam sequence (generalized Fibonacci sequence) $\{V_{n}(V_{0},V_{1};r,s)%
\}_{n\geq 0}$ (or shortly $\{V_{n}\}_{n\geq 0}$) is defined as follows: 
\begin{equation}
V_{n}=rV_{n-1}+sV_{n-2},\text{ \ \ \ \ }V_{0}=c_{0},V_{1}=c_{1},\text{ \ }%
n\geq 2  \label{equati:nbvvbuysdfacxz}
\end{equation}%
where $V_{0},V_{1}\ $are arbitrary reel or complex numbers and $r,s$ are
real numbers. The sequence $\{V_{n}\}_{n\geq 0}$ can be extended to negative
subscripts by defining%
\begin{equation*}
V_{-n}=-\frac{r}{s}V_{-(n-1)}+\frac{1}{s}V_{-(n-2)}
\end{equation*}%
for $n=1,2,3,...$ when $s\neq 0.$ Therefore, recurrence (\ref%
{equati:nbvvbuysdfacxz}) holds for all integer $n.$ See Table 5 for a few
members of Horadam sequences.

\ \ \ \ \ 

Table 5 A few members of Horadam sequences.

\begin{tabular}{cccc}
\hline
$\text{Sequences (Numbers)}$ &  &  & $\text{Notation}$ \\ \hline
$\text{Fibonacci}$ &  &  & $\{F_{n}\}=\{V_{n}(0,1;1,1)\}$ \\ 
$\text{Lucas}$ &  &  & $\{L_{n}\}=\{V_{n}(2,1;1,1)\}$ \\ 
$\text{Pell}$ &  &  & $\{P_{n}\}=\{V_{n}(0,1;2,1)\}$ \\ 
$\text{Pell-Lucas}$ &  &  & $\{Q_{n}\}=\{V_{n}(2,2;2,1)\}$ \\ 
$\text{second order Jacobsthal}$ &  &  & $\{J_{n}\}=\{V_{n}(0,1;1,2)\}$ \\ 
$\text{second order Jacobsthal-Lucas}$ &  &  & $\{j_{n}\}=\{V_{n}(2,1;1,2)\}$
\\ \hline
\end{tabular}

\ \ \ \ 

The first few values of the sequences with non-negative indices are shown
below (see Table 6).

\ \ 

Table 6 A few values of Horadam sequences with non-negative and negative
indices

\begin{tabular}{cccccccccccccccccc}
\hline
$n$ & $-8$ & $-7$ & $-6$ & $-5$ & $-4$ & $-3$ & $-2$ & $-1$ & $0$ & $1$ & $2$
& $3$ & $4$ & $5$ & $6$ & $7$ & $8$ \\ \hline
$F_{n}$ & $-21$ & $13$ & $-8$ & $5$ & $-3$ & $2$ & $-1$ & $1$ & $0$ & $1$ & $%
1$ & $2$ & $3$ & $5$ & $8$ & $13$ & $21$ \\ 
$L_{n}$ & $47$ & $-29$ & $18$ & $-11$ & $7$ & $-4$ & $3$ & $-1$ & $2$ & $1$
& $3$ & $4$ & $7$ & $11$ & $18$ & $29$ & $47$ \\ 
$P_{n}$ & $-408$ & $169$ & $-70$ & $29$ & $-12$ & $5$ & $-2$ & $1$ & $0$ & $%
1 $ & $2$ & $5$ & $12$ & $29$ & $70$ & $169$ & $408$ \\ 
$Q_{n}$ & $1154$ & $-478$ & $198$ & $-82$ & $34$ & $-14$ & $6$ & $-2$ & $2$
& $2$ & $6$ & $14$ & $34$ & $82$ & $198$ & $478$ & $1154$ \\ 
$J_{n}$ & $-\frac{85}{256}$ & $\frac{43}{128}$ & $-\frac{21}{64}$ & $\frac{11%
}{32}$ & $-\frac{5}{16}$ & $\frac{3}{8}$ & $-\frac{1}{4}$ & $\frac{1}{2}$ & $%
0$ & $1$ & $1$ & $3$ & $5$ & $11$ & $21$ & $43$ & $85$ \\ 
$j_{n}$ & $\frac{257}{256}$ & $-\frac{127}{128}$ & $\frac{65}{64}$ & $-\frac{%
31}{32}$ & $\frac{17}{16}$ & $-\frac{7}{8}$ & $\frac{5}{4}$ & $-\frac{1}{2}$
& $2$ & $1$ & $5$ & $7$ & $17$ & $31$ & $65$ & $127$ & $257$ \\ \hline
\end{tabular}

\ \ \ \ \ \ \ \ 

The generalized Tribonacci sequence $\{V_{n}(V_{0},V_{1},V_{2};r,s,t)\}_{n%
\geq 0}$ (or shortly $\{V_{n}\}_{n\geq 0}$) is defined as follows: 
\begin{equation}
V_{n}=rV_{n-1}+sV_{n-2}+tV_{n-3},\text{ \ \ \ \ }%
V_{0}=c_{0},V_{1}=c_{1},V_{2}=c_{2},\text{ \ }n\geq 3
\label{equation:hgfdstysopdcfxg}
\end{equation}%
where $V_{0},V_{1},V_{2}\ $are arbitrary reel or complexs numbers and $r,s,t$
are real numbers. The sequence $\{V_{n}\}_{n\geq 0}$ can be extended to
negative subscripts by defining%
\begin{equation*}
V_{-n}=-\frac{s}{t}V_{-(n-1)}-\frac{r}{t}V_{-(n-2)}+\frac{1}{t}V_{-(n-3)}
\end{equation*}%
for $n=1,2,3,...$ when $t\neq 0.$ Therefore, recurrence (\ref%
{equation:hgfdstysopdcfxg}) holds for all integer $n.$

In literature, for example, the following names and notations (see Table 7)
are used for the special case of $r,s,t$ and initial values.

\ \ \ \ \ 

Table 7 A few members of generalized Tribonacci sequences.

\begin{tabular}{cccc}
\hline
$\text{Sequences (Numbers)}$ &  &  & $\text{Notation}$ \\ \hline
$\text{Tribonacci}$ &  &  & $\{T_{n}\}=\{V_{n}(0,1,1;1,1,1)\}$ \\ 
$\text{Tribonacci-Lucas}$ &  &  & $\{K_{n}\}=\{V_{n}(3,1,3;1,1,1)\}$ \\ 
$\text{Padovan (Cordonnier)}$ &  &  & $\{P_{n}\}=\{V_{n}(1,1,1;0,1,1)\}$ \\ 
$\text{Pell-Padovan}$ &  &  & $\{R_{n}\}=\{V_{n}(1,1,1;0,2,1)\}$ \\ 
$\text{Jacobsthal-Padovan}$ &  &  & $\{JP_{n}\}=\{V_{n}(1,1,1;0,1,2)\}$ \\ 
$\text{Perrin}$ &  &  & $\{Q_{n}\}=\{V_{n}(3,0,2;0,1,1)\}$ \\ 
$\text{Pell-Perrin}$ &  &  & $\{pQ_{n}\}=\{V_{n}(3,0,2;0,2,1)\}$ \\ 
$\text{Jacobsthal-Perrin}$ &  &  & $\{JQ_{n}\}=\{V_{n}(3,0,2;0,1,2)\}$ \\ 
$\text{Padovan-Perrin}$ &  &  & $\{S_{n}\}=\{V_{n}(0,0,1;0,1,1)\}$ \\ 
$\text{Narayana}$ &  &  & $\{N_{n}\}=\{V_{n}(0,1,1;1,0,1)\}$ \\ 
$\text{third order Jacobsthal}$ &  &  & $\{J_{n}\}=\{V_{n}(0,1,1;1,1,2)\}$
\\ 
$\text{third order Jacobsthal-Lucas}$ &  &  & $\{j_{n}\}=%
\{V_{n}(2,1,5;1,1,2)\}$ \\ \hline
\end{tabular}

\ \ 

The first few values of the sequences with non-negative and negative indices
are shown below (see Table 8).

\ 

Table 8 A few values of generalized Tribonacci sequences.

\begin{tabular}{cccccccccccccccccc}
\hline
$n$ & $-8$ & $-7$ & $-6$ & $-5$ & $-4$ & $-3$ & $-2$ & $-1$ & $0$ & $1$ & $2$
& $3$ & $4$ & $5$ & $6$ & $7$ & $8$ \\ \hline
$T_{n}$ & $4$ & $1$ & $-3$ & $2$ & $0$ & $-1$ & $1$ & $0$ & $0$ & $1$ & $1$
& $2$ & $4$ & $7$ & $13$ & $24$ & $44$ \\ 
$K_{n}$ & $3$ & $-15$ & $11$ & $-1$ & $-5$ & $5$ & $-1$ & $-1$ & $3$ & $1$ & 
$3$ & $7$ & $11$ & $21$ & $39$ & $71$ & $131$ \\ 
$P_{n}$ & $0$ & $1$ & $-1$ & $1$ & $0$ & $0$ & $1$ & $0$ & $1$ & $1$ & $1$ & 
$2$ & $2$ & $3$ & $4$ & $5$ & $7$ \\ 
$R_{n}$ & $67$ & $-41$ & $25$ & $-15$ & $9$ & $-5$ & $3$ & $-1$ & $1$ & $1$
& $1$ & $3$ & $3$ & $7$ & $9$ & $17$ & $25$ \\ 
$JP_{n}$ & $\frac{23}{128}$ & $-\frac{3}{64}$ & $-\frac{1}{32}$ & $\frac{5}{%
16}$ & $-\frac{1}{8}$ & $\frac{1}{4}$ & $\frac{1}{2}$ & $0$ & $1$ & $1$ & $1$
& $3$ & $3$ & $5$ & $9$ & $11$ & $19$ \\ 
$Q_{n}$ & $5$ & $-1$ & $-2$ & $4$ & $-3$ & $2$ & $1$ & $-1$ & $3$ & $0$ & $2$
& $3$ & $2$ & $5$ & $5$ & $7$ & $10$ \\ 
$pQ_{n}$ & $156$ & $-96$ & $59$ & $-36$ & $22$ & $-13$ & $8$ & $-4$ & $3$ & $%
0$ & $2$ & $3$ & $4$ & $8$ & $11$ & $20$ & $30$ \\ 
$JQ_{n}$ & $\frac{161}{256}$ & $-\frac{85}{128}$ & $\frac{25}{64}$ & $\frac{%
19}{32}$ & $-\frac{15}{16}$ & $\frac{11}{8}$ & $\frac{1}{4}$ & $-\frac{1}{2}$
& $3$ & $0$ & $2$ & $6$ & $2$ & $10$ & $14$ & $14$ & $34$ \\ 
$S_{n}$ & $1$ & $-2$ & $2$ & $-1$ & $0$ & $1$ & $-1$ & $1$ & $0$ & $0$ & $1$
& $0$ & $1$ & $1$ & $1$ & $2$ & $2$ \\ 
$N_{n}$ & $0$ & $-2$ & $1$ & $1$ & $-1$ & $0$ & $1$ & $0$ & $0$ & $1$ & $1$
& $1$ & $2$ & $3$ & $4$ & $6$ & $9$ \\ 
$J_{n}$ & $\frac{55}{128}$ & $-\frac{9}{64}$ & $-\frac{9}{32}$ & $\frac{7}{16%
}$ & $-\frac{1}{8}$ & $-\frac{1}{4}$ & $\frac{1}{2}$ & $0$ & $0$ & $1$ & $1$
& $2$ & $5$ & $9$ & $18$ & $37$ & $73$ \\ 
$j_{n}$ & $-\frac{41}{32}$ & $\frac{7}{16}$ & $\frac{7}{8}$ & $-\frac{5}{4}$
& $\frac{1}{2}$ & $1$ & $-1$ & $1$ & $2$ & $1$ & $5$ & $10$ & $17$ & $37$ & $%
74$ & $145$ & $293$ \\ \hline
\end{tabular}%
\ \ 

\ \ \ \ 

The generalized Tetranacci sequence $%
\{V_{n}(V_{0},V_{1},V_{2},V_{3};r,s,t,u)\}_{n\geq 0}$ (or shortly $%
\{V_{n}\}_{n\geq 0}$) is defined as follows: 
\begin{equation}
V_{n}=rV_{n-1}+sV_{n-2}+tV_{n-3}+uV_{n-4},\text{ \ \ \ \ }%
V_{0}=c_{0},V_{1}=c_{1},V_{2}=c_{2},\text{ },V_{3}=c_{3},\text{\ }n\geq 4
\label{equati:dgvcxdtysdczx}
\end{equation}%
where $V_{0},V_{1},V_{2},V_{3}\ $are arbitrary reel or complex numbers and $%
r,s,t,u$ are real numbers. The sequence $\{V_{n}\}_{n\geq 0}$ can be
extended to negative subscripts by defining%
\begin{equation*}
V_{-n}=-\frac{t}{u}V_{-(n-1)}-\frac{s}{u}V_{-(n-2)}-\frac{r}{u}V_{-(n-3)}+%
\frac{1}{u}V_{-(n-4)}
\end{equation*}%
for $n=1,2,3,...$ when $u\neq 0.$ Therefore, recurrence (\ref%
{equati:dgvcxdtysdczx}) holds for all integer $n.$

In literature, for example, the following names and notations (see Table 9)
are used for the special case of $r,s,t,u$ and initial values.

\ \ \ \ \ 

Table 9 A few members of generalized Tetranacci sequences.

\begin{tabular}{cccc}
\hline
$\text{Sequences (Numbers)}$ &  &  & $\text{Notation}$ \\ \hline
$\text{T}$etranacci &  &  & $\{M_{n}\}=\{V_{n}(0,1,1,2;1,1,1,1)\}$ \\ 
Tetranacci$\text{-Lucas}$ &  &  & $\{R_{n}\}=\{V_{n}(4,1,3,7;1,1,1,1)\}$ \\ 
fourth$\text{ order Jacobsthal}$ &  &  & $\{J_{n}\}=\{V_{n}(0,1,1,1;1,1,1,2)%
\}$ \\ 
f$\text{ourth order Jacobsthal-Lucas}$ &  &  & $\{j_{n}\}=%
\{V_{n}(2,1,5,10;1,1,1,2)\}$ \\ \hline
\end{tabular}

\ \ 

The first few values of the sequences with non-negative and negative indices
are shown below (see Table 10).

\ 

Table 10 A few values of generalized Tetranacci sequences.

\begin{tabular}{ccccccccccccccccccc}
\hline
$n$ & $-8$ & $-7$ & $-6$ & $-5$ & $-4$ & $-3$ & $-2$ & $-1$ & $0$ & $1$ & $2$
& $3$ & $4$ & $5$ & $6$ & $7$ & $8$ & $9$ \\ \hline
$M_{n}$ & $-3$ & $2$ & $0$ & $0$ & $-1$ & $1$ & $0$ & $0$ & $0$ & $1$ & $1$
& $2$ & $4$ & $8$ & $15$ & $29$ & $56$ & $108$ \\ 
$R_{n}$ & $15$ & $-1$ & $-1$ & $-6$ & $7$ & $-1$ & $-1$ & $-1$ & $4$ & $1$ & 
$3$ & $7$ & $15$ & $26$ & $51$ & $99$ & $191$ & $367$ \\ 
$J_{n}$ & $-\frac{51}{256}$ & $\frac{77}{128}$ & $\frac{13}{64}$ & $-\frac{19%
}{32}$ & $-\frac{3}{16}$ & $\frac{5}{8}$ & $\frac{1}{4}$ & $-\frac{1}{2}$ & $%
0$ & $1$ & $1$ & $1$ & $3$ & $7$ & $13$ & $25$ & $51$ & $103$ \\ 
$j_{n}$ & $\frac{103}{128}$ & $-\frac{89}{64}$ & $\frac{7}{32}$ & $\frac{7}{%
16}$ & $\frac{7}{8}$ & $-\frac{5}{4}$ & $\frac{1}{2}$ & $1$ & $2$ & $1$ & $5$
& $10$ & $20$ & $37$ & $77$ & $154$ & $308$ & $613$ \\ \hline
\end{tabular}

\ \ \ \ 

The generalized Pentanacci sequence $%
\{V_{n}(V_{0},V_{1},V_{2},V_{3},V_{4};r,s,t,u,v)\}_{n\geq 0}$ (or shortly $%
\{V_{n}\}_{n\geq 0}$) is defined as follows: 
\begin{equation}
V_{n}=rV_{n-1}+sV_{n-2}+tV_{n-3}+uV_{n-4}+vV_{n-5},\text{\ }%
V_{0}=c_{0},V_{1}=c_{1},V_{2}=c_{2},\text{ }V_{3}=c_{3},V_{4}=c_{4},\text{\ }%
n\geq 5  \label{equati:gfdsdfvvxcvztysa}
\end{equation}%
where $V_{0},V_{1},V_{2},V_{3},V_{4}\ $are arbitrary reel or complex numbers
and $r,s,t,u$ are real numbers. The sequence $\{V_{n}\}_{n\geq 0}$ can be
extended to negative subscripts by defining%
\begin{equation*}
V_{-n}=-\frac{u}{v}V_{-n+1}-\frac{t}{v}V_{-n+2}-\frac{s}{v}V_{-n+3}-\frac{r}{%
v}V_{-n+4}+\frac{1}{v}V_{-n+5}
\end{equation*}%
for $n=1,2,3,...$ when $u\neq 0.$ Therefore, recurrence (\ref%
{equati:gfdsdfvvxcvztysa}) holds for all integer $n.$

In literature, for example, the following names and notations (see Table 11)
are used for the special case of $r,s,t,u,v$ and initial values.

\ \ \ \ \ 

Table 11 A few members of generalized Pentanacci sequences.

\begin{tabular}{cccc}
\hline
$\text{Sequences (Numbers)}$ &  &  & $\text{Notation}$ \\ \hline
Pentanacci &  &  & $\{P_{n}\}=\{V_{n}(0,1,1,2,4;1,1,1,1,1)\}$ \\ 
Pentanacci$\text{-Lucas}$ &  &  & $\{Q_{n}\}=\{V_{n}(5,1,3,7,15;1,1,1,1,1)\}$
\\ 
fifth$\text{ order Jacobsthal}$ &  &  & $\{J_{n}\}=%
\{V_{n}(0,1,1,1,1;1,1,1,1,2)\}$ \\ 
fifth$\text{ order Jacobsthal-Lucas}$ &  &  & $\{j_{n}\}=%
\{V_{n}(2,1,5,10,20;1,1,1,1,2)\}$ \\ \hline
\end{tabular}

\ \ 

The first few values of the sequences with non-negative and negative indices
are shown below (see Table 12).

\ 

Table 12 A few values of generalized Pentanacci sequences.

\begin{tabular}{ccccccccccccccccccc}
\hline
$n$ & $-8$ & $-7$ & $-6$ & $-5$ & $-4$ & $-3$ & $-2$ & $-1$ & $0$ & $1$ & $2$
& $3$ & $4$ & $5$ & $6$ & $7$ & $8$ & $9$ \\ \hline
$P_{n}$ & $0$ & $0$ & $0$ & $-1$ & $1$ & $0$ & $0$ & $0$ & $0$ & $1$ & $1$ & 
$2$ & $4$ & $8$ & $16$ & $31$ & $61$ & $120$ \\ 
$Q_{n}$ & $-1$ & $-1$ & $-7$ & $9$ & $-1$ & $-1$ & $-1$ & $-1$ & $5$ & $1$ & 
$3$ & $7$ & $15$ & $31$ & $57$ & $113$ & $223$ & $439$ \\ 
$J_{n}$ & $\frac{31}{64}$ & $-\frac{1}{32}$ & $-\frac{17}{16}$ & $-\frac{1}{8%
}$ & $\frac{3}{4}$ & $\frac{1}{2}$ & $0$ & $-1$ & $0$ & $1$ & $1$ & $1$ & $1$
& $4$ & $9$ & $17$ & $33$ & $65$ \\ 
$j_{n}$ & $\frac{13}{128}$ & $\frac{13}{64}$ & $\frac{13}{32}$ & $\frac{13}{%
16}$ & $-\frac{11}{8}$ & $\frac{1}{4}$ & $\frac{1}{2}$ & $1$ & $2$ & $1$ & $%
5 $ & $10$ & $20$ & $40$ & $77$ & $157$ & $314$ & $628$ \\ \hline
\end{tabular}

\ 

\section{Particular Cases of Main Result}

There is a well-known Simson Identity (formula) for Fibonacci sequence $%
\{F_{n}\}$, namely,%
\begin{equation*}
F_{n+1}F_{n-1}-F_{n}^{2}=(-1)^{n}
\end{equation*}%
which was derived first by R. Simson in 1753 and it is now called as Cassini
Identity (formula) as well. This can be written in the form%
\begin{equation*}
\left\vert 
\begin{array}{cc}
F_{n+1} & F_{n} \\ 
F_{n} & F_{n-1}%
\end{array}%
\right\vert =(-1)^{n}.
\end{equation*}

A search of the literature turns up that there are many identities including
Simson (Cassini), Catalan, d'Ocagne, Melham, Tagiuri, Gelin-Cesaro, Gould
identities, see for example, [\ref{cooper2015some},\ref{fairgrieve2005pro},%
\ref{hendel2017pro},\ref{koshy2015gelin},\ref{lang2013afib},\ref%
{lang2013bfib},\ref{melham2003afibon},\ref{melham2011onpro}].

Next, we consider generalized Horadam numbers $V_{n}=rV_{n-1}+sV_{n-2}$ with 
$2$ initial terms $V_{0}=c_{0},$ $V_{1}=c_{1}$ and present a formula for
those numbers.

\begin{theorem}[Simson Formula of Horadam Numbers]
\label{theorem:hgbcdszaos}For all integers $n$ we have%
\begin{equation}
\left\vert 
\begin{array}{cc}
V_{n+1} & V_{n} \\ 
V_{n} & V_{n-1}%
\end{array}%
\right\vert =(-1)^{n}s^{n}\left\vert 
\begin{array}{cc}
V_{1} & V_{0} \\ 
V_{0} & V_{-1}%
\end{array}%
\right\vert .  \label{equation:gfdfghopskmnfs}
\end{equation}
\end{theorem}

Proof. We proof by induction on $n.$ Firstly, we prove the formula (\ref%
{equation:gfdfghopskmnfs}) for $n\geq 0.$ For $n=0,$ it is obvious that the
formula is true. Now, we assume that the formula (\ref%
{equation:gfdfghopskmnfs}) is true for $n=k,$ that is 
\begin{equation*}
\left\vert 
\begin{array}{cc}
V_{k+1} & V_{k} \\ 
V_{k} & V_{k-1}%
\end{array}%
\right\vert =(-1)^{k}s^{k}\left\vert 
\begin{array}{cc}
V_{1} & V_{0} \\ 
V_{0} & V_{-1}%
\end{array}%
\right\vert .
\end{equation*}%
Then by induction hypothesis, we obtain%
\begin{eqnarray*}
\left\vert 
\begin{array}{cc}
V_{k+2} & V_{k+1} \\ 
V_{k+1} & V_{k}%
\end{array}%
\right\vert &=&\left\vert 
\begin{array}{cc}
rV_{k+1}+sV_{k} & V_{k+1} \\ 
rV_{k}+sV_{k-1} & V_{k}%
\end{array}%
\right\vert =r\left\vert 
\begin{array}{cc}
V_{k+1} & V_{k+1} \\ 
V_{k} & V_{k}%
\end{array}%
\right\vert +s\left\vert 
\begin{array}{cc}
V_{k} & V_{k+1} \\ 
V_{k-1} & V_{k}%
\end{array}%
\right\vert \\
&=&-s\left\vert 
\begin{array}{cc}
V_{k+1} & V_{k} \\ 
V_{k} & V_{k-1}%
\end{array}%
\right\vert =-s\left( (-1)^{k}s^{k}\left\vert 
\begin{array}{cc}
V_{1} & V_{0} \\ 
V_{0} & V_{-1}%
\end{array}%
\right\vert \right) \\
&=&(-1)^{k+1}s^{k+1}\left\vert 
\begin{array}{cc}
V_{1} & V_{0} \\ 
V_{0} & V_{-1}%
\end{array}%
\right\vert .
\end{eqnarray*}%
i.e., the formula (\ref{equation:gfdfghopskmnfs}) is true for $n=k+1.$ Thus,
(\ref{equation:gfdfghopskmnfs}) hold for all integers $n\geq 1.$

Now we consider the formula (\ref{equation:gfdfghopskmnfs}) for $n\leq -1.$
Take $h=-n$ so that $h\geq 1$. So we need to prove by induction that for $%
h\geq 1$ we have%
\begin{equation}
\left\vert 
\begin{array}{cc}
V_{-h+1} & V_{-h} \\ 
V_{-h} & V_{-h-1}%
\end{array}%
\right\vert =(-1)^{-h}s^{-h}\left\vert 
\begin{array}{cc}
V_{1} & V_{0} \\ 
V_{0} & V_{-1}%
\end{array}%
\right\vert .  \label{equation:rtgfcvasdfgb}
\end{equation}%
For $h=1,$ the formula is true because%
\begin{eqnarray*}
\left\vert 
\begin{array}{cc}
V_{0} & V_{-1} \\ 
V_{-1} & V_{-2}%
\end{array}%
\right\vert &=&-\left\vert 
\begin{array}{cc}
V_{-1} & V_{0} \\ 
V_{-2} & V_{-1}%
\end{array}%
\right\vert =-\left\vert 
\begin{array}{cc}
-\frac{r}{s}V_{0}+\frac{1}{s}V_{1} & V_{0} \\ 
-\frac{r}{s}V_{-1}+\frac{1}{s}V_{0} & V_{-1}%
\end{array}%
\right\vert \\
&=&-\left\vert 
\begin{array}{cc}
-\frac{r}{s}V_{0} & V_{0} \\ 
-\frac{r}{s}V_{-1} & V_{-1}%
\end{array}%
\right\vert -\left\vert 
\begin{array}{cc}
\frac{1}{s}V_{1} & V_{0} \\ 
\frac{1}{s}V_{0} & V_{-1}%
\end{array}%
\right\vert =-\frac{1}{s}\left\vert 
\begin{array}{cc}
V_{1} & V_{0} \\ 
V_{0} & V_{-1}%
\end{array}%
\right\vert .
\end{eqnarray*}%
Now, we assume that the formula (\ref{equation:rtgfcvasdfgb}) is true for $%
h=k,$ that is%
\begin{equation}
\left\vert 
\begin{array}{cc}
V_{-k+1} & V_{-k} \\ 
V_{-k} & V_{-k-1}%
\end{array}%
\right\vert =(-1)^{-k}s^{-k}\left\vert 
\begin{array}{cc}
V_{1} & V_{0} \\ 
V_{0} & V_{-1}%
\end{array}%
\right\vert .  \label{equation:sdfghbvcxrtyuds}
\end{equation}%
Then by induction hypothesis (\ref{equation:sdfghbvcxrtyuds}), we obtain%
\begin{eqnarray*}
\left\vert 
\begin{array}{cc}
V_{-(k+1)+1} & V_{-(k+1)} \\ 
V_{-(k+1)} & V_{-(k+1)-1}%
\end{array}%
\right\vert &=&\left\vert 
\begin{array}{cc}
V_{-k} & V_{-k-1} \\ 
V_{-k-1} & V_{-k-2}%
\end{array}%
\right\vert =-\left\vert 
\begin{array}{cc}
V_{-k-1} & V_{-k} \\ 
V_{-k-2} & V_{-k-1}%
\end{array}%
\right\vert =-\left\vert 
\begin{array}{cc}
-\frac{r}{s}V_{-k}+\frac{1}{s}V_{-k+1} & V_{-k} \\ 
-\frac{r}{s}V_{-k-1}+\frac{1}{s}V_{-k} & V_{-k-1}%
\end{array}%
\right\vert \\
&=&-\left\vert 
\begin{array}{cc}
-\frac{r}{s}V_{-k} & V_{-k} \\ 
-\frac{r}{s}V_{-k-1} & V_{-k-1}%
\end{array}%
\right\vert -\left\vert 
\begin{array}{cc}
\frac{1}{s}V_{-k+1} & V_{-k} \\ 
\frac{1}{s}V_{-k} & V_{-k-1}%
\end{array}%
\right\vert \\
&=&-\frac{1}{s}\left\vert 
\begin{array}{cc}
V_{-k+1} & V_{-k} \\ 
V_{-k} & V_{-k-1}%
\end{array}%
\right\vert =-\frac{1}{s}\left( (-1)^{-k}s^{-k}\left\vert 
\begin{array}{cc}
V_{1} & V_{0} \\ 
V_{0} & V_{-1}%
\end{array}%
\right\vert \right) \\
&=&(-1)^{-k+1}s^{-(k+1)}\left\vert 
\begin{array}{cc}
V_{1} & V_{0} \\ 
V_{0} & V_{-1}%
\end{array}%
\right\vert =(-1)^{-(k+1)}s^{-(k+1)}\left\vert 
\begin{array}{cc}
V_{1} & V_{0} \\ 
V_{0} & V_{-1}%
\end{array}%
\right\vert
\end{eqnarray*}%
so that the formula (\ref{equation:rtgfcvasdfgb}) is true for $h=k+1.$ Thus,
(\ref{equation:rtgfcvasdfgb}) holds for all integers $h\geq 1$ and so (\ref%
{equation:gfdfghopskmnfs}) holds for all integers $n\leq -1.$ This completes
the proof.%
%TCIMACRO{\TeXButton{End Proof}{\endproof}}%
%BeginExpansion
\endproof%
%EndExpansion

\begin{remark}
Theorem \ref{theorem:hgbcdszaos} is given in Horadam [\ref%
{horadam1965generfunc}] (see also [\ref{horadam1965basicprope}]). In fact,
in [\ref{horadam1965generfunc}], Horadam gave a beautiful formula more
general case, namely Catalan Identity for Horadam numbers. We provide the
proof of Theorem \ref{theorem:hgbcdszaos} here because it pave the way the
method to prove the general case.
\end{remark}

\ \ \ We can write Theorem\ \ref{theorem:hgbcdszaos} as 
\begin{equation*}
f(n)=(-1)^{n}s^{n}f(0)
\end{equation*}%
where $f(n)=\left\vert 
\begin{array}{cc}
V_{n+1} & V_{n} \\ 
V_{n} & V_{n-1}%
\end{array}%
\right\vert $ and $f(0)=\left\vert 
\begin{array}{cc}
V_{1} & V_{0} \\ 
V_{0} & V_{-1}%
\end{array}%
\right\vert $. \ In the following Table 13, we present Simsons's formula of
particular\ Horadam sequences.

\ \ \ 

Table 13 Simsons's formula of some Horadam sequences

\begin{tabular}{ccccccc}
\hline
Sequence: $V_{n}$ &  & Simson Formula &  & Sequence: $V_{n}$ &  & Simson
Formula \\ \hline
$F_{n}$ &  & $f(n)=(-1)^{n}$ &  & $L_{n}$ &  & $f(n)=5(-1)^{n-1}$ \\ 
$P_{n}$ &  & $f(n)=(-1)^{n}$ &  & $Q_{n}$ &  & $f(n)=8(-1)^{n-1}$ \\ 
$J_{n}$ &  & $f(n)=(-1)^{n}2^{n-1}$ &  & $j_{n}$ &  & $%
f(n)=9(-1)^{n-1}2^{n-1}$ \\ \hline
\end{tabular}

\ \ \ \ 

Next we consider generalized Tribonacci numbers $%
V_{n}=rV_{n-1}+sV_{n-2}+tV_{n-3}$ with $3$ initial terms $V_{0}=c_{0},$ $%
V_{1}=c_{1},$ $V_{2}=c_{2}.$

\begin{theorem}[Simson Formula of Generalized Tribonacci Numbers]
\label{theorem:tribnmcgxftsyudx}For all integers $n$ we have%
\begin{equation}
\left\vert 
\begin{array}{ccc}
V_{n+2} & V_{n+1} & V_{n} \\ 
V_{n+1} & V_{n} & V_{n-1} \\ 
V_{n} & V_{n-1} & V_{n-2}%
\end{array}%
\right\vert =t^{n}\left\vert 
\begin{array}{ccc}
V_{2} & V_{1} & V_{0} \\ 
V_{1} & V_{0} & V_{-1} \\ 
V_{0} & V_{-1} & V_{-2}%
\end{array}%
\right\vert .  \label{equation:mehzehazocxb}
\end{equation}
\end{theorem}

Proof. We prove by induction on $n.$ Firstly, we prove the formula (\ref%
{equation:mehzehazocxb}) for $n\geq 0.$ For $n=0,$ it is obvious that the
formula is true. Now, we assume that the formula (\ref{equation:mehzehazocxb}%
) is true for $n=k,$ that is 
\begin{equation*}
\left\vert 
\begin{array}{ccc}
V_{k+2} & V_{k+1} & V_{k} \\ 
V_{k+1} & V_{k} & V_{k-1} \\ 
V_{k} & V_{k-1} & V_{k-2}%
\end{array}%
\right\vert =t^{k}\left\vert 
\begin{array}{ccc}
V_{2} & V_{1} & V_{0} \\ 
V_{1} & V_{0} & V_{-1} \\ 
V_{0} & V_{-1} & V_{-2}%
\end{array}%
\right\vert .
\end{equation*}%
Then by induction hypothesis, we obtain%
\begin{eqnarray*}
&&\left\vert 
\begin{array}{ccc}
V_{k+3} & V_{k+2} & V_{k+1} \\ 
V_{k+2} & V_{k+1} & V_{k} \\ 
V_{k+1} & V_{k} & V_{k-1}%
\end{array}%
\right\vert \\
&=&\left\vert 
\begin{array}{ccc}
rV_{k+2}+sV_{k+1}+tV_{k} & V_{k+2} & V_{k+1} \\ 
rV_{k+1}+sV_{k}+tV_{k-1} & V_{k+1} & V_{k} \\ 
rV_{k}+sV_{k-1}+tV_{k-2} & V_{k} & V_{k-1}%
\end{array}%
\right\vert \\
&=&\left\vert 
\begin{array}{ccc}
rV_{k+2} & V_{k+2} & V_{k+1} \\ 
rV_{k+1} & V_{k+1} & V_{k} \\ 
rV_{k} & V_{k} & V_{k-1}%
\end{array}%
\right\vert +\left\vert 
\begin{array}{ccc}
sV_{k+1} & V_{k+2} & V_{k+1} \\ 
sV_{k} & V_{k+1} & V_{k} \\ 
sV_{k-1} & V_{k} & V_{k-1}%
\end{array}%
\right\vert +\left\vert 
\begin{array}{ccc}
tV_{k} & V_{k+2} & V_{k+1} \\ 
tV_{k-1} & V_{k+1} & V_{k} \\ 
tV_{k-2} & V_{k} & V_{k-1}%
\end{array}%
\right\vert \\
&=&t\left\vert 
\begin{array}{ccc}
V_{k} & V_{k+2} & V_{k+1} \\ 
V_{k-1} & V_{k+1} & V_{k} \\ 
V_{k-2} & V_{k} & V_{k-1}%
\end{array}%
\right\vert =t\left\vert 
\begin{array}{ccc}
V_{k+2} & V_{k+1} & V_{k} \\ 
V_{k+1} & V_{k} & V_{k-1} \\ 
V_{k} & V_{k-1} & V_{k-2}%
\end{array}%
\right\vert \\
&=&t\left( t^{k}\left\vert 
\begin{array}{ccc}
V_{2} & V_{1} & V_{0} \\ 
V_{1} & V_{0} & V_{-1} \\ 
V_{0} & V_{-1} & V_{-2}%
\end{array}%
\right\vert \right) =t^{k+1}\left\vert 
\begin{array}{ccc}
V_{2} & V_{1} & V_{0} \\ 
V_{1} & V_{0} & V_{-1} \\ 
V_{0} & V_{-1} & V_{-2}%
\end{array}%
\right\vert
\end{eqnarray*}%
i.e., the formula (\ref{equation:mehzehazocxb}) is true for $n=k+1.$ Thus, (%
\ref{equation:mehzehazocxb}) holds for all integers $n\geq 1.$

Now we consider the formula (\ref{equation:mehzehazocxb}) for $n\leq -1.$
Take $h=-n$ so that $h\geq 1$. So we need to prove by induction that for $%
h\geq 1$ \ 
\begin{equation}
\left\vert 
\begin{array}{ccc}
V_{-h+2} & V_{-h+1} & V_{-h} \\ 
V_{-h+1} & V_{-h} & V_{-h-1} \\ 
V_{-h} & V_{-h-1} & V_{-h-2}%
\end{array}%
\right\vert =t^{-h}\left\vert 
\begin{array}{ccc}
V_{2} & V_{1} & V_{0} \\ 
V_{1} & V_{0} & V_{-1} \\ 
V_{0} & V_{-1} & V_{-2}%
\end{array}%
\right\vert .  \label{equation:ybvszkardomnub}
\end{equation}%
For $h=1,$ the formula is true because%
\begin{eqnarray*}
\left\vert 
\begin{array}{ccc}
V_{1} & V_{0} & V_{-1} \\ 
V_{0} & V_{-1} & V_{-2} \\ 
V_{-1} & V_{-2} & V_{-3}%
\end{array}%
\right\vert &=&\left\vert 
\begin{array}{ccc}
V_{-1} & V_{1} & V_{0} \\ 
V_{-2} & V_{0} & V_{-1} \\ 
V_{-3} & V_{-1} & V_{-2}%
\end{array}%
\right\vert =\left\vert 
\begin{array}{ccc}
-\frac{s}{t}V_{0}-\frac{r}{t}V_{1}+\allowbreak \frac{1}{t}V_{2} & V_{1} & 
V_{0} \\ 
-\frac{s}{t}V_{-1}-\frac{r}{t}V_{0}+\allowbreak \frac{1}{t}V_{1} & V_{0} & 
V_{-1} \\ 
-\frac{s}{t}V_{-2}-\frac{r}{t}V_{-1}+\allowbreak \frac{1}{t}V_{0} & V_{-1} & 
V_{-2}%
\end{array}%
\right\vert \\
&=&\left\vert 
\begin{array}{ccc}
-\frac{s}{t}V_{0} & V_{1} & V_{0} \\ 
-\frac{s}{t}V_{-1} & V_{0} & V_{-1} \\ 
-\frac{s}{t}V_{-2} & V_{-1} & V_{-2}%
\end{array}%
\right\vert +\left\vert 
\begin{array}{ccc}
-\frac{r}{t}V_{1} & V_{1} & V_{0} \\ 
-\frac{r}{t}V_{0} & V_{0} & V_{-1} \\ 
-\frac{r}{t}V_{-1} & V_{-1} & V_{-2}%
\end{array}%
\right\vert +\left\vert 
\begin{array}{ccc}
\allowbreak \frac{1}{t}V_{2} & V_{1} & V_{0} \\ 
\allowbreak \frac{1}{t}V_{1} & V_{0} & V_{-1} \\ 
\allowbreak \frac{1}{t}V_{0} & V_{-1} & V_{-2}%
\end{array}%
\right\vert \\
&=&\allowbreak \frac{1}{t}\left\vert 
\begin{array}{ccc}
V_{2} & V_{1} & V_{0} \\ 
V_{1} & V_{0} & V_{-1} \\ 
V_{0} & V_{-1} & V_{-2}%
\end{array}%
\right\vert .
\end{eqnarray*}%
Now, we assume that the formula (\ref{equation:ybvszkardomnub}) is true for $%
h=k,$ that is 
\begin{equation}
\left\vert 
\begin{array}{ccc}
V_{-k+2} & V_{-k+1} & V_{-k} \\ 
V_{-k+1} & V_{-k} & V_{-k-1} \\ 
V_{-k} & V_{-k-1} & V_{-k-2}%
\end{array}%
\right\vert =t^{-k}\left\vert 
\begin{array}{ccc}
V_{2} & V_{1} & V_{0} \\ 
V_{1} & V_{0} & V_{-1} \\ 
V_{0} & V_{-1} & V_{-2}%
\end{array}%
\right\vert .  \label{equation:ghesewqazxtyf}
\end{equation}%
Then by induction hypothesis, we obtain%
\begin{eqnarray*}
&&\left\vert 
\begin{array}{ccc}
V_{-(k+1)+2} & V_{-(k+1)+1} & V_{-(k+1)} \\ 
V_{-(k+1)+1} & V_{-(k+1)} & V_{-(k+1)-1} \\ 
V_{-(k+1)} & V_{-(k+1)-1} & V_{-(k+1)-2}%
\end{array}%
\right\vert =\left\vert 
\begin{array}{ccc}
V_{-k+1} & V_{-k} & V_{-k-1} \\ 
V_{-k} & V_{-k-1} & V_{-k-2} \\ 
V_{-k-1} & V_{-k-2} & V_{-k-3}%
\end{array}%
\right\vert \\
&=&\left\vert 
\begin{array}{ccc}
V_{-k-1} & V_{-k+1} & V_{-k} \\ 
V_{-k-2} & V_{-k} & V_{-k-1} \\ 
V_{-k-3} & V_{-k-1} & V_{-k-2}%
\end{array}%
\right\vert =\left\vert 
\begin{array}{ccc}
-\frac{s}{t}V_{-k}-\frac{r}{t}V_{-k+1}+\allowbreak \frac{1}{t}V_{-k+2} & 
V_{-k+1} & V_{-k} \\ 
-\frac{s}{t}V_{-k-1}-\frac{r}{t}V_{-k}+\allowbreak \frac{1}{t}V_{-k+1} & 
V_{-k} & V_{-k-1} \\ 
-\frac{s}{t}V_{-k-2}-\frac{r}{t}V_{-k-1}+\allowbreak \frac{1}{t}V_{-k} & 
V_{-k-1} & V_{-k-2}%
\end{array}%
\right\vert \\
&=&\left\vert 
\begin{array}{ccc}
-\frac{s}{t}V_{-k} & V_{-k+1} & V_{-k} \\ 
-\frac{s}{t}V_{-k-1} & V_{-k} & V_{-k-1} \\ 
-\frac{s}{t}V_{-k-2} & V_{-k-1} & V_{-k-2}%
\end{array}%
\right\vert +\left\vert 
\begin{array}{ccc}
-\frac{r}{t}V_{-k+1} & V_{-k+1} & V_{-k} \\ 
-\frac{r}{t}V_{-k} & V_{-k} & V_{-k-1} \\ 
-\frac{r}{t}V_{-k-1} & V_{-k-1} & V_{-k-2}%
\end{array}%
\right\vert +\left\vert 
\begin{array}{ccc}
\allowbreak \frac{1}{t}V_{-k+2} & V_{-k+1} & V_{-k} \\ 
\allowbreak \frac{1}{t}V_{-k+1} & V_{-k} & V_{-k-1} \\ 
\allowbreak \frac{1}{t}V_{-k} & V_{-k-1} & V_{-k-2}%
\end{array}%
\right\vert \\
&=&\frac{1}{t}\left\vert 
\begin{array}{ccc}
\allowbreak V_{-k+2} & V_{-k+1} & V_{-k} \\ 
\allowbreak V_{-k+1} & V_{-k} & V_{-k-1} \\ 
V_{-k} & V_{-k-1} & V_{-k-2}%
\end{array}%
\right\vert =\frac{1}{t}t^{-k}\left\vert 
\begin{array}{ccc}
V_{2} & V_{1} & V_{0} \\ 
V_{1} & V_{0} & V_{-1} \\ 
V_{0} & V_{-1} & V_{-2}%
\end{array}%
\right\vert =t^{-(k+1)}\left\vert 
\begin{array}{ccc}
V_{2} & V_{1} & V_{0} \\ 
V_{1} & V_{0} & V_{-1} \\ 
V_{0} & V_{-1} & V_{-2}%
\end{array}%
\right\vert
\end{eqnarray*}%
i.e., the formula (\ref{equation:ybvszkardomnub}) is true for $h=k+1.$ Thus,
(\ref{equation:ybvszkardomnub}) holds for all integers $h\geq 1$ and so (\ref%
{equation:mehzehazocxb}) holds for all integers $n\leq -1.$ This completes
the proof. 
%TCIMACRO{\TeXButton{End Proof}{\endproof}}%
%BeginExpansion
\endproof%
%EndExpansion

\ \ \ We can write Theorem\ \ref{theorem:tribnmcgxftsyudx} as 
\begin{equation*}
f(n)=t^{n}f(0)
\end{equation*}%
where $f(n)=\left\vert 
\begin{array}{ccc}
V_{n+2} & V_{n+1} & V_{n} \\ 
V_{n+1} & V_{n} & V_{n-1} \\ 
V_{n} & V_{n-1} & V_{n-2}%
\end{array}%
\right\vert $ and $f(0)=\left\vert 
\begin{array}{ccc}
V_{2} & V_{1} & V_{0} \\ 
V_{1} & V_{0} & V_{-1} \\ 
V_{0} & V_{-1} & V_{-2}%
\end{array}%
\right\vert $.

\ \ \ 

In the following Table 14, we present Simsons's formula of particular\
generalized Tribonacci sequences.

\ \ 

Table 14 Simsons's formula of some generalized Tribonacci sequences

\begin{tabular}{ccccccc}
\hline
Sequence: $V_{n}$ &  & Simson Formula &  & Sequence: $V_{n}$ &  & Simson
Formula \\ \hline
$T_{n}$ &  & $f(n)=-1$ &  & $K_{n}$ &  & $f(n)=-44$ \\ 
$P_{n}$ &  & $f(n)=-1$ &  & $R_{n}$ &  & $f(n)=-4$ \\ 
$JP_{n}$ &  & $f(n)=-2^{n}$ &  & $Q_{n}$ &  & $f(n)=-23$ \\ 
$pQ_{n}$ &  & $f(n)=-11$ &  & $JQ_{n}$ &  & $f(n)=-13\times 2^{n+1}$ \\ 
$S_{n}$ &  & $f(n)=-1$ &  & $N_{n}$ &  & $f(n)=-1$ \\ 
$J_{n}$ &  & $f(n)=-2^{n-1}$ &  & $j_{n}$ &  & $f(n)=-9\times 2^{n+1}$ \\ 
\hline
\end{tabular}

\ \ \ 

Next we consider generalized Tetranacci numbers $%
V_{n}=rV_{n-1}+sV_{n-2}+tV_{n-3}+uV_{n-4}$ with $4$ initial terms%
\begin{equation*}
V_{0}=c_{0},\text{ }V_{1}=c_{1},\text{ }V_{2}=c_{2},V_{3}=c_{3}.\text{ }
\end{equation*}

\begin{theorem}[Simson Formula of Generalized Tetranacci Numbers]
\label{theorem:sdrtqopkmnvbxc}For all integers $n$ we have%
\begin{equation}
\left\vert 
\begin{array}{cccc}
V_{n+3} & V_{n+2} & V_{n+1} & V_{n} \\ 
V_{n+2} & V_{n+1} & V_{n} & V_{n-1} \\ 
V_{n+1} & V_{n} & V_{n-1} & V_{n-2} \\ 
V_{n} & V_{n-1} & V_{n-2} & V_{n-3}%
\end{array}%
\right\vert =(-1)^{n}u^{n}\left\vert 
\begin{array}{cccc}
V_{3} & V_{2} & V_{1} & V_{0} \\ 
V_{2} & V_{1} & V_{0} & V_{-1} \\ 
V_{1} & V_{0} & V_{-1} & V_{-2} \\ 
V_{0} & V_{-1} & V_{-2} & V_{-3}%
\end{array}%
\right\vert .  \label{equation:hghgfgvbcuopdfcx}
\end{equation}
\end{theorem}

\textit{Proof}. The proof can be given exactly as the proof of Theorem \ref%
{theorem:hgbcdszaos}, so we omit it. 
%TCIMACRO{\TeXButton{End Proof}{\endproof}}%
%BeginExpansion
\endproof%
%EndExpansion

\ \ \ We can write Theorem\ \ref{theorem:sdrtqopkmnvbxc} as 
\begin{equation*}
f(n)=(-1)^{n}u^{n}f(0)
\end{equation*}%
where $f(n)=\left\vert 
\begin{array}{cccc}
V_{n+3} & V_{n+2} & V_{n+1} & V_{n} \\ 
V_{n+2} & V_{n+1} & V_{n} & V_{n-1} \\ 
V_{n+1} & V_{n} & V_{n-1} & V_{n-2} \\ 
V_{n} & V_{n-1} & V_{n-2} & V_{n-3}%
\end{array}%
\right\vert $ and $f(0)=\left\vert 
\begin{array}{cccc}
V_{3} & V_{2} & V_{1} & V_{0} \\ 
V_{2} & V_{1} & V_{0} & V_{-1} \\ 
V_{1} & V_{0} & V_{-1} & V_{-2} \\ 
V_{0} & V_{-1} & V_{-2} & V_{-3}%
\end{array}%
\right\vert $.

\ \ 

In the following Table 15, we present Simsons's formula of particular\
generalized Tetranacci sequences.

\ \ \ 

Table 15 Simsons's formula of some generalized Tetranacci sequences

\begin{tabular}{ccc}
\hline
Sequence: $V_{n}$ &  & Simson Formula \\ \hline
$M_{n}$ &  & $f(n)=(-1)^{n-1}$ \\ 
$R_{n}$ &  & $f(n)=563(-1)^{n}$ \\ 
$J_{n}$ &  & $f(n)=0$ \\ 
$j_{n}$ &  & $f(n)=(-1)^{n}2^{n-2}3^{5}$ \\ \hline
\end{tabular}

\ \ 

Next we consider generalized Pentanacci numbers $%
V_{n}=rV_{n-1}+sV_{n-2}+tV_{n-3}+uV_{n-4}+vV_{n-5}$ with $5$ initial terms $%
V_{0}=c_{0},$ $V_{1}=c_{1},$ $V_{2}=c_{2},V_{3}=c_{3},V_{4}=c_{4}.$

\begin{theorem}[Simson Formula of Generalized Pentanacci Numbers]
\label{theorem:dfgcvzrtqweadsz}For all integers $n$ we have%
\begin{equation}
\left\vert 
\begin{array}{ccccc}
V_{n+4} & V_{n+3} & V_{n+2} & V_{n+1} & V_{n} \\ 
V_{n+3} & V_{n+2} & V_{n+1} & V_{n} & V_{n-1} \\ 
V_{n+2} & V_{n+1} & V_{n} & V_{n-1} & V_{n-2} \\ 
V_{n+1} & V_{n} & V_{n-1} & V_{n-2} & V_{n-3} \\ 
V_{n} & V_{n-1} & V_{n-2} & V_{n-3} & V_{n-4}%
\end{array}%
\right\vert =v^{n}\left\vert 
\begin{array}{ccccc}
V_{4} & V_{3} & V_{2} & V_{1} & V_{0} \\ 
V_{3} & V_{2} & V_{1} & V_{0} & V_{-1} \\ 
V_{2} & V_{1} & V_{0} & V_{-1} & V_{-2} \\ 
V_{1} & V_{0} & V_{-1} & V_{-2} & V_{-3} \\ 
V_{0} & V_{-1} & V_{-2} & V_{-3} & V_{-4}%
\end{array}%
\allowbreak \right\vert .  \label{equation:haduohbesdaz}
\end{equation}
\end{theorem}

\textit{Proof}. The proof can be given exactly as the proof of Theorem \ref%
{theorem:tribnmcgxftsyudx}, so we omit it. 
%TCIMACRO{\TeXButton{End Proof}{\endproof}}%
%BeginExpansion
\endproof%
%EndExpansion

\ \ \ We can write Theorem\ \ref{theorem:dfgcvzrtqweadsz} as 
\begin{equation*}
f(n)=v^{n}f(0)
\end{equation*}%
where $f(n)=\left\vert 
\begin{array}{ccccc}
V_{n+4} & V_{n+3} & V_{n+2} & V_{n+1} & V_{n} \\ 
V_{n+3} & V_{n+2} & V_{n+1} & V_{n} & V_{n-1} \\ 
V_{n+2} & V_{n+1} & V_{n} & V_{n-1} & V_{n-2} \\ 
V_{n+1} & V_{n} & V_{n-1} & V_{n-2} & V_{n-3} \\ 
V_{n} & V_{n-1} & V_{n-2} & V_{n-3} & V_{n-4}%
\end{array}%
\right\vert $ and $f(0)=\left\vert 
\begin{array}{ccccc}
V_{4} & V_{3} & V_{2} & V_{1} & V_{0} \\ 
V_{3} & V_{2} & V_{1} & V_{0} & V_{-1} \\ 
V_{2} & V_{1} & V_{0} & V_{-1} & V_{-2} \\ 
V_{1} & V_{0} & V_{-1} & V_{-2} & V_{-3} \\ 
V_{0} & V_{-1} & V_{-2} & V_{-3} & V_{-4}%
\end{array}%
\right\vert $.

\ \ 

In the following Table 16, we present Simsons's formula of particular\
generalized Pentanacci sequences.

\ \ \ \ 

Table 16 Simsons's formula of some generalized Pentanacci sequences

\begin{tabular}{ccc}
\hline
Sequence: $V_{n}$ &  & Simson Formula \\ \hline
$P_{n}$ &  & $f(n)=1$ \\ 
$Q_{n}$ &  & $f(n)=9584$ \\ 
$J_{n}$ &  & $f(n)=2^{n-2}\times 11$ \\ 
$j_{n}$ &  & $f(n)=2^{n-3}\times 3^{4}\times 19$ \\ \hline
\end{tabular}

\ \ 

\section{Main Result}

Now we consider the $m$-order linear recurrence relation%
\begin{equation*}
V_{n}=%
\sum_{i=1}^{m}r_{i}V_{n-i}=r_{1}V_{n-1}+r_{2}V_{n-2}+r_{3}V_{n-3}+...+r_{m}V_{n-m}.
\end{equation*}%
For $m\geq 2$,\ we define $f$ by%
\begin{equation*}
f(n)=\left\vert 
\begin{array}{ccccccc}
V_{n+m-1} & V_{n+m-2} & V_{n+m-3} & \cdots & V_{n+2} & V_{n+1} & V_{n} \\ 
V_{n+m-2} & V_{n+m-3} & V_{n+m-4} & \cdots & V_{n+1} & V_{n} & V_{n-1} \\ 
V_{n+m-3} & V_{n+m-4} & V_{n+m-5} & \cdots & V_{n} & V_{n-1} & V_{n-2} \\ 
\vdots & \vdots & \vdots & \vdots & \vdots & \vdots & \vdots \\ 
V_{n+2} & V_{n+1} & V_{n} & \cdots & V_{n-m+5} & V_{n-m+4} & V_{n-m+3} \\ 
V_{n+1} & V_{n} & V_{n-1} & \cdots & V_{n-m+4} & V_{n-m+3} & V_{n-m+2} \\ 
V_{n} & V_{n-1} & V_{n-2} & \cdots & V_{n-m+3} & V_{n-m+2} & V_{n-m+1}%
\end{array}%
\right\vert .
\end{equation*}%
Note that 
\begin{equation*}
f(0)=\left\vert 
\begin{array}{ccccccc}
V_{m-1} & V_{m-2} & V_{m-3} & \cdots & V_{2} & V_{1} & V_{0} \\ 
V_{m-2} & V_{m-3} & V_{m-4} & \cdots & V_{1} & V_{0} & V_{-1} \\ 
V_{m-3} & V_{m-4} & V_{m-5} & \cdots & V_{0} & V_{-1} & V_{-2} \\ 
\vdots & \vdots & \vdots & \vdots & \vdots & \vdots & \vdots \\ 
V_{2} & V_{1} & V_{0} & \cdots & V_{m+5} & V_{-m+4} & V_{-m+3} \\ 
V_{1} & V_{0} & V_{-1} & \cdots & V_{m+4} & V_{-m+3} & V_{-m+2} \\ 
V_{0} & V_{-1} & V_{-2} & \cdots & V_{m+3} & V_{-m+2} & V_{-m+1}%
\end{array}%
\right\vert .
\end{equation*}

Motivated by the cases $m=2,3,4,5,$ we are ready to present our main result
for the arbitrary $m.$

\begin{theorem}[Simson Formula of Generalized $m$-step Fibonacci Numbers]
\label{theorem:hgfdsyunbvgh}Let $m\geq 2$. Then for all integers $n$ we have 
\begin{equation}
f(n)=y(n)r_{m}^{n}f(0)  \label{equation:tysbgvxczdfnbc}
\end{equation}%
where%
\begin{equation*}
y(n)=\left\{ 
\begin{array}{ccc}
1 & , & m\text{ odd} \\ 
(-1)^{n} & , & m\text{ even}%
\end{array}%
\right. .
\end{equation*}
\end{theorem}

\textit{Proof.} We prove the theorem by induction for $n\geq 0,$ the proof
of the case $n\leq -1$ being similar. As in the proof of the cases $%
m=2,3,4,5 $ we need to consider $m$ separately as odd and even. We provide
the proof of the even cases. For $n=0,$ it is obvious that the formula is
true. Now, we assume that the formula (\ref{equation:tysbgvxczdfnbc}) is
true for $n=k.$ Then we will complete the inductive step $n=k+1$ as follows:
Note that%
\begin{equation*}
f(k+1)=\left\vert 
\begin{array}{ccccccc}
V_{k+m} & V_{k+m-1} & V_{k+m-2} & \cdots & V_{k+3} & V_{k+2} & V_{k+1} \\ 
V_{k+m-1} & V_{k+m-2} & V_{k+m-3} & \cdots & V_{k+2} & V_{k+1} & V_{k} \\ 
V_{k+m-2} & V_{k+m-3} & V_{k+m-4} & \cdots & V_{k+1} & V_{k} & V_{k} \\ 
\vdots & \vdots & \vdots & \vdots & \vdots & \vdots & \vdots \\ 
V_{k+3} & V_{k+2} & V_{k+1} & \cdots & V_{k-m+6} & V_{k-m+5} & V_{k-m+4} \\ 
V_{k+2} & V_{k+1} & V_{k} & \cdots & V_{k-m+5} & V_{k-m+4} & V_{k-m+3} \\ 
V_{k+1} & V_{k} & V_{k-1} & \cdots & V_{k-m+4} & V_{k-m+3} & V_{k-m+2}%
\end{array}%
\right\vert .
\end{equation*}%
Using the recurrence relations%
\begin{eqnarray*}
V_{k+m} &=&r_{1}V_{k+m-1}+r_{2}V_{k+m-2}+r_{3}V_{k+m-3}+...+r_{m}V_{k} \\
V_{k+m-1} &=&r_{1}V_{k+m-2}+r_{2}V_{k+m-3}+r_{3}V_{k+m-4}+...+r_{m}V_{k-1} \\
V_{k+m-2} &=&r_{1}V_{k+m-3}+r_{2}V_{k+m-4}+r_{3}V_{k+m-5}+...+r_{m}V_{k-2} \\
&&\vdots \\
V_{k+3} &=&r_{1}V_{k+2}+r_{2}V_{k+1}+r_{3}V_{k}+...+r_{m}V_{k-m+3} \\
V_{k+2} &=&r_{1}V_{k+1}+r_{2}V_{k}+r_{3}V_{k-1}+...+r_{m}V_{k-m+2} \\
V_{k+1} &=&r_{1}V_{k}+r_{2}V_{k-1}+r_{3}V_{k+1-3}+...+r_{m}V_{k-m+1}
\end{eqnarray*}%
in the $1^{st}$ column of the determinant $f(k+1)$ and expanding $1^{st}$
column as $m-1$ additions and then after rearranging the determinant, we
obtain%
\begin{eqnarray*}
\ \text{\ }\ f(k+1) &=&r_{m}\left\vert 
\begin{array}{ccccccc}
V_{k} & V_{k+m-1} & V_{k+m-2} & \cdots & V_{k+3} & V_{k+2} & V_{k+1} \\ 
V_{k-1} & V_{k+m-2} & V_{k+m-3} & \cdots & V_{k+2} & V_{k+1} & V_{k} \\ 
V_{k-2} & V_{k+m-3} & V_{k+m-4} & \cdots & V_{k+1} & V_{k} & V_{k} \\ 
\vdots & \vdots & \vdots & \vdots & \vdots & \vdots & \vdots \\ 
V_{k-m+3} & V_{k+2} & V_{k+1} & \cdots & V_{k-m+6} & V_{k-m+5} & V_{k-m+4}
\\ 
V_{k-m+2} & V_{k+1} & V_{k} & \cdots & V_{k-m+5} & V_{k-m+4} & V_{k-m+3} \\ 
V_{k+1-m} & V_{k} & V_{k-1} & \cdots & V_{k-m+4} & V_{k-m+3} & V_{k-m+2}%
\end{array}%
\right\vert \\
&=&-r_{m}\left\vert 
\begin{array}{ccccccc}
V_{k+m-1} & V_{k+m-2} & V_{k+m-3} & \cdots & V_{k+2} & V_{k+1} & V_{k} \\ 
V_{k+m-2} & V_{k+m-3} & V_{k+m-4} & \cdots & V_{k+1} & V_{k} & V_{k-1} \\ 
V_{k+m-3} & V_{k+m-4} & V_{k+m-5} & \cdots & V_{k} & V_{k-1} & V_{k-2} \\ 
\vdots & \vdots & \vdots & \vdots & \vdots & \vdots & \vdots \\ 
V_{k+2} & V_{k+1} & V_{k} & \cdots & V_{k-m+5} & V_{k-m+4} & V_{k-m+3} \\ 
V_{k+1} & V_{k} & V_{k-1} & \cdots & V_{k-m+4} & V_{k-m+3} & V_{k-m+2} \\ 
V_{k} & V_{k-1} & V_{k-2} & \cdots & V_{k-m+3} & V_{k-m+2} & V_{k-m+1}%
\end{array}%
\right\vert \\
&=&-r_{m}((-1)^{k}r_{m}^{k}f(0))=(-1)^{k+1}r_{m}^{k+1}f(0).
\end{eqnarray*}%
This completes the inductive step and the proof of the theorem. 
%TCIMACRO{\TeXButton{End Proof}{\endproof}}%
%BeginExpansion
\endproof%
%EndExpansion

\begin{remark}
Of course, this paper could be shorthened. To calculate Simson Identity we
needed sequences and the values of the elements of those sequences. But a
search of the literature shows that it is not easy to find sequences of
altogether the case $m=2,3,4,5\ $of the generalized $m$-step Fibonacci
numbers in a single reference. So, as much as presenting new results, we
wanted to fill this gap as well by giving the sequences and the values of
their elements as tables.
\end{remark}

\end{document}